\documentclass[reqno]{amsart}
\usepackage{amssymb,hyperref}

\theoremstyle{plain}
\newtheorem{thm}{Theorem}
\newtheorem{lem}{Lemma}
\newtheorem{cor}{Corollary}

\theoremstyle{definition}
\newtheorem{ex}{Example}

\renewcommand{\Re}{\mathrm{Re}}

\title
[Some classes of order $\alpha$ for second-order differential inequalities]
{Some classes of order $\alpha$ \\
for second-order differential inequalities}

\author{Hitoshi Shiraishi}
\address{Hitoshi Shiraishi \newline
Department of Mathematics \newline
Kinki University \newline
Higashi-Osaka, Osaka 577-8502, Japan}
\email{step\_625@hotmail.com}

\author{Kazuo Kuroki}
\address{Kazuo Kuroki \newline
Department of Mathematics \newline
Kinki University \newline
Higashi-Osaka, Osaka 577-8502, Japan}
\email{freedom@sakai.zaq.ne.jp}

\author{Shigeyoshi Owa}
\address{Shigeyoshi Owa \newline
Department of Mathematics \newline
Kinki University \newline
Higashi-Osaka, Osaka 577-8502, Japan}
\email{shige21@ican.zaq.ne.jp}


\date{}

\begin{document}

\begin{abstract}
For analytic functions $f(z)$ in the open unit disk $\mathbb{U}$ with $f(0)=f'(0)-1=0$,
S. S. Miller and P. T. Mocanu (Integral Transform. Spec. Funct. {\bf 19}(2008)) have considered some sufficient problems for starlikeness.
The object of the present paper is to discuss some sufficient problems for $f(z)$ to be in some classes of order $\alpha$.
\end{abstract}

\begin{flushleft}
This paper was published in the journal: \\
Electron. J. Math. Anal. Appl. {\bf 1} (2013), No. 2, 149--155. \\
{\footnotesize \url{http://ejmaa.6te.net/Vol2_Papers/Vol.%201(2)%20July%202013,%20pp.%20149-155..pdf}}
\end{flushleft}
\hrule

\

\

\maketitle

\section{Introduction}

\

Let $\mathcal{A}_n$ denote the class of functions
$$
f(z)=z+a_{n+1}z^{n+1}+a_{n+2}z^{n+2}+ \ldots
\qquad(n=1,2,3,\ldots)
$$
that are analytic in the open unit disk $\mathbb{U}=\{z \in \mathbb{C}:|z|<1\}$
and $\mathcal{A}=\mathcal{A}_1$.
We denote by $\mathcal{S}$ the subclass of $\mathcal{A}_n$ consisting of univalent functions $f(z)$ in $\mathbb{U}$.

Let $\mathcal{S^{*}}(\alpha)$ be defined by
$$
\mathcal{S^{*}}(\alpha)
=\left\{f(z)\in\mathcal{A}_n:
\Re\left(\frac{zf'(z)}{f(z)}\right)>\alpha,\ 
0\leqq {}^{\exists} \alpha<1 \right\}.
$$
We denote by $\mathcal{S}^{*} = \mathcal{S}^{*}(0)$.
Also,
let $\mathcal{C}(\alpha)$ be
$$
\mathcal{C}({\alpha})
=\left\{f(z)\in\mathcal{A}_n:
\Re(f'(z))>\alpha,\ 
0\leqq {}^{\exists} \alpha<1 \right\}.
$$
We also denote by $\mathcal{C}=\mathcal{C}(0)$.
Also,
let $\mathcal{K}(\alpha)$ be defined by
$$
\mathcal{K}(\alpha)
=\left\{f(z)\in\mathcal{A}_n:
\Re\left(1+\frac{zf''(z)}{f'(z)}\right)>\alpha,\ 
0\leqq {}^{\exists} \alpha<1 \right\}.
$$
We denote by $\mathcal{K}=\mathcal{K}(0)$.
From the definitions for $\mathcal{S^{*}}(\alpha)$ and $\mathcal{K}(\alpha)$,
we know that $f(z) \in \mathcal{K}(\alpha)$ if and only if $zf'(z) \in \mathcal{S}^{*}(\alpha)$.

\

The basic tool in proving our results is the following lemma due to Jack \cite{m1ref1}
(also, due to Miller and Mocanu \cite{m1ref2}).

\

\begin{lem} \label{jack} \quad
Let the function $w(z)$ defined by
$$
w(z)=a_nz^n+a_{n+1}z^{n+1}+a_{n+2}z^{n+2}+ \ldots
\qquad(n=1,2,3,\ldots)
$$
be analytic in $\mathbb{U}$ with $w(0)=0$.
If $\left|w(z)\right|$ attains its maximum value on the circle $\left|z\right|=r$ at a point $z_{0}\in\mathbb{U}$,
then there exists a real number $k \geqq n$ such that
$$
\frac{z_{0}w'(z_{0})}{w(z_{0})}=k
$$
and
$$
\Re\left(\frac{z_0w''(z_0)}{w'(z_0)}\right)+1\geqq k.
$$
\end{lem}

\

\section{Main results}

\

Applying Lemma \ref{jack},
we derive the following lemma.

\

\begin{lem} \label{d1lem1.1} \quad
If $f(z)\in\mathcal{A}_n$ satisfies
$$
\left|zf''(z)-\beta\left(f'(z)-\frac{f(z)}{z}\right)\right|
< \rho|n+1-\beta|
\qquad(z\in\mathbb{U})
$$
for some real $\rho>0$
and some complex $\beta$
with $\Re(\beta)<n+1$,
then
$$
\left|f'(z)-\frac{f(z)}{z}\right|<\rho
\qquad (z\in\mathbb{U}).
$$
\end{lem}

\

\begin{proof} \quad
Let us define $w(z)$ by
\begin{align}
w(z)
&= f'(z)-\frac{f(z)}{z} \label{d1lem1.1eq1}\\
&= n a_{n+1}z^{n} + (n+1)a_{n+2}z^{n+1} + \ldots
\qquad(z\in\mathbb{U}). \nonumber
\end{align}
Then, clearly,
$w(z)$ is analytic in $\mathbb{U}$ and $w(0)=0$.
Differentiating both sides in (\ref{d1lem1.1eq1}),
we obtain
$$
zf''(z)=zw'(z)+w(z)
\qquad(z\in\mathbb{U}),
$$
and therefore,
\begin{align*}
\left|zf''(z)-\beta\left(f'(z)-\frac{f(z)}{z}\right)\right|
&= |zw'(z)+(1-\beta)w(z)|\\
&= |w(z)| \left|\frac{zw'(z)}{w(z)}+1-\beta \right|\\
&< \rho|n+1-\beta|
\qquad(z\in\mathbb{U}).
\end{align*}
If there exists a point $z_{0} \in \mathbb{U}$ such that
$$
\max_{\left| z \right| \leqq \left| z_{0} \right|} \left| w(z) \right|
= \left| w(z_{0}) \right|
= \rho,
$$
then Lemma \ref{jack} gives us that
$w(z_{0})=\rho e^{i \theta}$ and $z_{0}w'(z_{0})=kw(z_{0})$ ($k \geqq n$).
Thus we have
\begin{align*}
\left|z_0f''(z_0)-\beta\left(f'(z_0)-\frac{f(z_0)}{z_0}\right)\right|
&= |w(z_0)| \left|\frac{z_0w'(z_0)}{w(z_0)}+1-\beta \right|\\
&= \rho|k+1-\beta|\\
&\geqq \rho|n+1-\beta|.
\end{align*}
This contradicts our condition in the lemma.
Therefore,
there is no $z_{0} \in \mathbb{U}$ such that $\left|w(z_{0})\right| = \rho$.
This means that $\left|w(z)\right| < \rho$ for all $z \in \mathbb{U}$,
that is, that
$$
\left| f'(z)-\frac{f(z)}{z} \right| < \rho
\qquad(z\in\mathbb{U}).
$$
\end{proof}

\

Also applying Lemma \ref{jack},
we have

\

\begin{lem} \label{d1lem1.2} \quad
If $f(z)\in\mathcal{A}_n$ satisfies
$$
\left|zf''(z) -\beta\left(f'(z) -\frac{f(z)}{z}\right)\right|
<\rho n|n+1-\beta|
\qquad(z\in\mathbb{U})
$$
for some real $\rho>0$
and some complex $\beta$
with $\Re(\beta)<n+1$,
then
$$
\left|\frac{f(z)}{z}-1\right|
<\rho
\qquad(z\in\mathbb{U}).
$$
\end{lem}

\

\begin{proof} \quad
Let us define the function $w(z)$ by
\begin{align*}
w(z) &= \frac{f(z)}{z}-1\\
&= a_{n+1}z^{n} + a_{n+2}z^{n+1} + \ldots
\qquad(z\in\mathbb{U}).
\end{align*}
Clearly,
$w(z)$ is analytic in $\mathbb{U}$ and $w(0)=0$.
We want to prove that $|w(z)|<\rho$ in $\mathbb{U}$.
Since
$$
zf''(z)=z^2w''(z)+2zw'(z)
\qquad(z\in\mathbb{U}),
$$
we see that
\begin{align*}
\left|zf''(z) -\beta\left(f'(z) -\frac{f(z)}{z}\right)\right|
&= |z^2w''(z) +(2-\beta)zw'(z)|\\
&< \rho n|n+1-\beta|
\qquad(z\in\mathbb{U}).
\end{align*}
If there exists a point $z_{0} \in \mathbb{U}$ such that
$$
\max_{\left| z \right| \leqq \left| z_{0} \right|} \left| w(z) \right|
= \left| w(z_{0}) \right|
= \rho,
$$
then Lemma \ref{jack} gives us that
$w(z_{0})=\rho e^{i \theta}$, $z_{0}w'(z_{0})=kw(z_{0})$ ($k \geqq n$) and
$$
\Re\left(\dfrac{z_0w''(z_0)}{w'(z_0)}\right)+1
\geqq k.
$$
Thus we have
\begin{align*}
\left|z_0f''(z_0) -\beta\left(f'(z_0) -\frac{f(z_0)}{z_0}\right)\right|
&= |z_0^2w''(z_0) +(2-\beta)z_0w'(z_0)|\\
&= |z_0w'(z_0)|\left|\frac{z_0w''(z_0)}{w'(z_0)}+2-\beta\right|\\
&= \rho k\left|\frac{z_0w''(z_0)}{w'(z_0)}+2-\beta\right|\\
&\geqq \rho k\left|\Re\left(\frac{z_0w''(z_0)}{w'(z_0)}\right)+2-\beta\right|\\
&\geqq \rho k|k+1-\beta|\\
&\geqq \rho n|n+1-\beta|.
\end{align*}
This contradicts the condition in the lemma.
Therefore, there is no $z_{0}\in\mathbb{U}$ such that $|w(z_{0})|=\rho$.
This means that $|w(z)| < \rho$ for all $z\in\mathbb{U}$.
\end{proof}

\

From Lemma \ref{d1lem1.1} and Lemma \ref{d1lem1.2},
we drive the following results for $\mathcal{S^{*}}(\alpha)$.

\

\begin{thm} \label{d1thm1} \quad
If $f(z)\in\mathcal{A}_n$ satisfies
$$
\left|zf''(z) -\beta\left(f'(z) -\frac{f(z)}{z}\right)\right|
<\frac{(1-\alpha)n|n+1-\beta|}{n+1-\alpha}
\qquad( z\in\mathbb{U})
$$
for some real $0\leqq\alpha<1$
and some complex $\beta$
with $\Re(\beta)<n+1$,
then
$$
\left|\frac{zf'(z)}{f(z)}-1\right|
< 1-\alpha
\qquad(z\in\mathbb{U}),
$$
so that $f(z)\in\mathcal{S}^{*}(\alpha)$.
\end{thm}

\

\begin{proof} \quad
From Lemma \ref{d1lem1.1} and Lemma \ref{d1lem1.2},
we have
\begin{equation}
\left| f'(z)-\frac{f(z)}{z} \right|
< \frac{n(1-\alpha)}{n+1-\alpha}
\qquad(z\in\mathbb{U}) \label{d1thm1eq1}
\end{equation}
and
\begin{equation}
\left| \frac{f(z)}{z}-1 \right|
< \frac{1-\alpha}{n+1-\alpha}
\qquad(z\in\mathbb{U}). \label{d1thm1eq2}
\end{equation}
From (\ref{d1thm1eq1}) and (\ref{d1thm1eq2}),
\begin{align*}
\frac{n(1-\alpha)}{n+1-\alpha}
&> \left|f'(z)-\frac{f(z)}{z}\right|\\
&= \left|\frac{f(z)}{z}\right|\left|\frac{zf'(z)}{f(z)}-1\right|\\
&> \left(1 -\frac{1-\alpha}{n+1-\alpha}\right)\left|\frac{zf'(z)}{f(z)}-1\right|\\
&= \frac{n}{n+1-\alpha}\left|\frac{zf'(z)}{f(z)}-1\right|
\qquad (z\in\mathbb{U}).
\end{align*}
So,
we can get
$$
\frac{n}{n+1-\alpha}\left|\frac{zf'(z)}{f(z)}-1\right|
< \frac{n(1-\alpha)}{n+1-\alpha}
\qquad(z\in\mathbb{U}).
$$
which completes the proof of the theorem.
\end{proof}

\

When we put $f(z)$ by $zf'(z)$ in Theorem \ref{d1thm1},
we have

\

\begin{cor} \label{d1cor1} \quad
If $f(z)\in\mathcal{A}_n$ satisfies
$$
|z^2f'''(z) +(2-\beta)zf''(z)|
<\frac{(1-\alpha)n|n+1-\beta|}{n+1-\alpha}
\qquad(z\in\mathbb{U})
$$
for some real $0\leqq\alpha<1$
and some complex $\beta$
with $\Re(\beta)<n+1$,
then
$$
\left|\left(1+\frac{zf''(z)}{f'(z)}\right)-1\right|
< 1-\alpha
\qquad(z\in\mathbb{U}),
$$
so that $f(z)\in\mathcal{K}(\alpha)$.
\end{cor}

\

\begin{ex} \label{d1ex1} \quad
For some real $0\leqq\alpha<1$
and some complex $\beta$
with $\Re(\beta)<n+1$,
we consider the function $f(z)$ given by
$$
f(z)=z+\frac{1-\alpha}{n+1-\alpha}z^{n+1}
\qquad(z\in\mathbb{U}).
$$
The function $f(z)$ satisfies Theorem \ref{d1thm1}.
\end{ex}

\

Next,
we consider $\mathcal{C}(\alpha)$.

\

\begin{thm} \label{d1thm2} \quad
If $f(z)\in\mathcal{A}_n$ satisfies
$$
|zf''(z)-\beta(f'(z)-1)|
<(1-\alpha)|n-\beta|
\qquad(z\in\mathbb{U})
$$
for some real $0\leqq\alpha<1$
and some complex $\beta$
with $\Re(\beta)<n$,
then
$$
|f'(z)-1|<1-\alpha
\qquad(z\in\mathbb{U}).
$$
This means that $f(z)\in\mathcal{C}(\alpha).$
\end{thm}

\

\begin{proof} \quad
Define $w(z)$ in $\mathbb{U}$ by
\begin{align}
w(z)
&= \frac{f'(z)-1}{1-\alpha} \label{d1thm2eq1}\\
&= \frac{(n+1)a_{n+1}}{1-\alpha}z^{n} + \frac{(n+2)a_{n+2}}{1-\alpha}z^{n+1} + \ldots
\qquad(z\in\mathbb{U}). \nonumber
\end{align}
Evidently,
$w(z)$ analytic in $\mathbb{U}$ and $w(0)=0$.
We want to prove $|w(z)|<1$.
Differentiating (\ref{d1thm2eq1}) and simplifiying,
we obtain
$$
zf''(z)=(1-\alpha)zw'(z)
\qquad (z\in\mathbb{U}).
$$
and,
hence
\begin{align*}
|zf''(z)-\beta(f'(z)-1)|
&= |(1-\alpha)zw'(z)-\beta(1-\alpha)w(z)|\\
&= (1-\alpha)|w(z)|\left|\frac{zw'(z)}{w(z)}-\beta\right|\\
&< (1-\alpha)|n-\beta|
\qquad(z\in\mathbb{U}).
\end{align*}
If there exists a point $z_{0} \in \mathbb{U}$ such that
$$
\max_{\left| z \right| \leqq \left| z_{0} \right|} \left| w(z) \right|
= \left| w(z_{0}) \right|
= 1,
$$
then Lemma \ref{jack} gives us that
$w(z_{0})= e^{i \theta}$ and $z_{0}w'(z_{0})=kw(z_{0})$ ($k \geqq n$).
Thus we have
\begin{align*}
|z_0f''(z_0)-\beta(f'(z_0)-1)|
&= (1-\alpha)|w(z_0)|\left|\frac{z_0w'(z_0)}{w(z_0)}-\beta\right|\\
&= (1-\alpha)|k-\beta|\\
&\geqq (1-\alpha)|n-\beta|.
\end{align*}
This contradicts our condition in the theorem.
Therefore,
there is no $z_0\in\mathbb{U}$ such that $w(z_0)=1$.
This means that $|w(z)|<1$ for all $z\in\mathbb{U}$.
\end{proof}

\

\begin{ex} \label{d1ex2} \quad
For some real $0\leqq\alpha<1$
and some complex $\beta$
with $\Re(\beta)<n$,
we take
$$
f(z)=z+\frac{1-\alpha}{n+1}z^{n+1}
\qquad(z\in\mathbb{U}).
$$
Then,
$f(z)$ satisfies Theorem \ref{d1thm2}.
\end{ex}

\

We get the following lemma from Lemma \ref{jack}.

\

\begin{lem} \label{d1lem3} \quad
If $f(z)\in\mathcal{A}_n$ satisfies
$$
|zf''(z)-\beta(f'(z)-1)|
<\rho|n-\beta|
\qquad(z\in\mathbb{U})
$$
for some real $\rho>0$
and some complex $\beta$
with $\Re(\beta)<n$,
then
$$
|f'(z)-1|<\rho
\qquad(z\in\mathbb{U}).
$$
\end{lem}

\

\begin{proof} \quad
Letting
\begin{align*}
w(z) &= f'(z)-1\\
&= (n+1)a_{n+1}z^{n} + (n+2)a_{n+2}z^{n+1} + \ldots
\qquad(z\in\mathbb{U}),
\end{align*}
we see that $w(z)$ is analytic in $\mathbb{U}$ and $w(0)=0$.
Noting that
$$
zf''(z)=zw'(z)
\qquad(z\in\mathbb{U}),
$$
we have
\begin{align*}
|zf''(z)-\beta(f'(z)-1)|
&= |zw'(z)-\beta w(z)|\\
&= |w(z)|\left|\frac{zw'(z)}{w(z)}-\beta\right|\\
&< \rho|n-\beta|
\qquad(z\in\mathbb{U}).
\end{align*}
If there exists a point $z_{0} \in \mathbb{U}$ such that
$$
\max_{\left| z \right| \leqq \left| z_{0} \right|} \left| w(z) \right|
= \left| w(z_{0}) \right|
= \rho,
$$
then Lemma \ref{jack} gives us that
$w(z_{0})=\rho e^{i \theta}$ and $z_{0}w'(z_{0})=kw(z_{0})$ ($k \geqq n$).
Thus we have
\begin{align*}
|z_0f''(z_0)-\beta(f'(z_0)-1)|
&= |w(z_0)|\left|\frac{z_0w'(z_0)}{w(z_0)}-\beta\right|\\
&= \rho|k-\beta|\\
&\geqq \rho|n-\beta|
\end{align*}
which contradicts our condition in the lemma.
Therefore,
there is no $z_0\in\mathbb{U}$ such that $|w(z_0)|=\rho$.
This means that $|w(z)|<\rho$ for all $z\in\mathbb{U}$.
\end{proof}

\

Using Lemma \ref{d1lem3},
we have next theorem.

\

\begin{thm} \label{d1thm3} \quad
If $f(z)\in\mathcal{A}_n$ satisfies
$$
|zf''(z)-\beta(f'(z)-1)|
<\alpha|n-\beta|
\qquad(z\in\mathbb{U})
$$
for some real $0<\alpha\leqq\dfrac{1}{2}$
and some complex $\beta$
with $\Re(\beta)<n$,
or
$$
|zf''(z)-\beta(f'(z)-1)|
<(1-\alpha)|n-\beta|
\qquad(z\in\mathbb{U})
$$
for some real $\dfrac{1}{2}\leqq\alpha<1$
and some complex $\beta$
with $\Re(\beta)<n$,
then
$$
\left|\frac{1}{f'(z)}-\frac{1}{2\alpha}\right|
<\frac{1}{2\alpha}
\qquad(z\in\mathbb{U}),
$$
which implies that $f(z)\in\mathcal{C}(\alpha)$.
\end{thm}

\

\begin{proof} \quad
We can get
\begin{equation}
|f'(z)-1|<\rho
\qquad(z\in\mathbb{U}).
\label{d1thm3eq1}
\end{equation}
for $0<\alpha \leqq \dfrac{1}{2}$ and $\rho=\alpha$,
or $\dfrac{1}{2} \leqq \alpha<1$ and $\rho=1-\alpha$
from Lemma \ref{d1lem3}.
Using (\ref{d1thm3eq1}),
we have
$$
|f'(z)-2\alpha|<S<|f'(z)|
$$
for $0<\alpha \leqq \dfrac{1}{2}$ and $S=1-\alpha$,
or $\dfrac{1}{2} \leqq \alpha<1$ and $S=\alpha$.
Thus we get
\begin{align*}
S\left| \frac{1}{f'(z)} -\frac{1}{2\alpha} \right|
&< |f'(z)| \left| \frac{1}{f'(z)}-\frac{1}{2\alpha} \right|\\
&= \left| 1-\frac{f'(z)}{2\alpha} \right|\\
&= \frac{1}{2\alpha} |f'(z)-2\alpha|\\
&< \frac{S}{2\alpha}
\qquad(z\in\mathbb{U}).
\end{align*}
So we obtain
$$
S\left| \frac{1}{f'(z)}-\frac{1}{2\alpha} \right|
< \frac{S}{2\alpha}
\qquad(z\in\mathbb{U}).
$$
\end{proof}

\

\begin{ex} \label{d1ex3} \quad
For some real $0<\alpha\leqq\dfrac{1}{2}$
and some complex $\beta$
with $\Re(\beta)<n$,
we consider the function $f(z)$ given by
$$
f(z)=z+\frac{\alpha}{n+1}z^{n+1}
\qquad(z\in\mathbb{U}).
$$
The function $f(z)$ satisfies Theorem \ref{d1thm3}.
\end{ex}

\

\end{document}